\documentclass[letterpaper, 10 pt, conference,fleqn]{ieeeconf}

\IEEEoverridecommandlockouts                           
\usepackage{times,latexsym,amssymb,amsmath,theorem}
\usepackage[hyperindex]{hyperref}
\usepackage{graphicx}
\usepackage{amssymb,amsmath}

\title{\LARGE \bf
Rational Observers of Rational Systems
}
\author{Jana N\v{e}mcov\'{a} and Mih\'{a}ly Petreczky and Jan H. van Schuppen
\thanks{This work was partially supported by the GA\v{C}R project 13-16764P.}
\thanks{J. N\v{e}mcov\'{a} is with the Department of Mathematics, University of Chemistry and Technology, Technick\'{a} 5, 166 28 Prague 6, Czech republic
        {\tt\small Jana.Nemcova@vscht.cz}}%
\thanks{M. Petreczky is with the Ecole Centrale Lille, 59651 Villeneuve d'Asq, Lille, France
        {\tt\small mihaly.petreczky@ec-lille.fr}}%
\thanks{J.H. van Schuppen is with the company Van Schuppen Control Research, 
          Gouden Leeuw 143, 1103 KB Amsterdam, The Netherlands
        {\tt\small  jan.h.van.schuppen@xs4all.nl}}%
}

\newtheorem{theorem}{Theorem}[section]

\newtheorem{definition}[theorem]{Definition}

\newtheorem{problem}[theorem]{Problem}
\newtheorem{procedure}[theorem]{Procedure}
\newtheorem{proposition}[theorem]{Proposition}

{\theorembodyfont{\rmfamily} \newtheorem{example}[theorem]{Example}}

\begin{document}

\newcommand{\sign}{\mbox{sign}}

\newcommand{\mathbbd}{\mbox{${\mathbb{D}}$}} 
\newcommand{\DS}{\mbox{${\mathbb{DS}}$}} 
\newcommand{\DSC}{\mbox{${\mathbb{DSC}}$}} 
\newcommand{\M}{\mbox{${\mathbb{M}}$}} 
\newcommand{\N}{\mbox{${\mathbb{N}}$}} 
\newcommand{\PP}{\mbox{${\mathbb{P}}$}} 
\newcommand{\Q}{\mbox{${\mathbb{Q}}$}} 
\newcommand{\W}{\mbox{${\mathbb{W}}$}} 
\newcommand{\C}{\mbox{$C$}} 
\newcommand{\R}{\mbox{$R$}} 
\newcommand{\T}{\mbox{$T$}} 
\newcommand{\Z}{\mbox{$Z$}}

\newcommand{\boundary}{\mbox{$\rm bndry$}}
\newcommand{\card}{\mbox{$\rm card$}}
\newcommand{\cisetsmin}{\mbox{$\rm CISets_{\min}$}}
\newcommand{\cisets}{\mbox{$\rm CISets$}}
\newcommand{\cix}{\mbox{$\rm CIX$}}
\newcommand{\colpwr}{\mbox{Colpwr}}
\newcommand{\closure}{\mbox{closure}}
\newcommand{\da}{\mbox{$\rm DA$}}
\newcommand{\domainofattraction}{\mbox{$DA$}}
\newcommand{\ds}{\mbox{$\rm DS$}}
\newcommand{\graph}{\mbox{Graph}}
\newcommand{\init}{\mbox{Init}}
\newcommand{\interior}{\mbox{int}}
\newcommand{\Mp}{\mbox{\bf M}}
\newcommand{\powerset}{\mbox{Pwrset}}
\newcommand{\pwr}{\mbox{Pwrset}}

\newcommand{\orderlexg}{\mbox{$>_{{\rm lex}}$}}
\newcommand{\orderlexl}{\mbox{$<_{{\rm lex}}$}}
\newcommand{\nm}{\mbox{$\mathbb{N}_m$}}
\newcommand{\relationconnect}{\mbox{${\rm R_{connect}}$}}
\newcommand{\zk}{\mbox{$\mathbb{Z}_k$}}
\newcommand{\zn}{\mbox{$\mathbb{Z}_n$}}
\newcommand{\zm}{\mbox{$\mathbb{Z}_m$}}
\newcommand{\zo}{\mbox{${\rm O}$}} 
\newcommand{\zr}{\mbox{$\mathbb{Z}_r$}}
\newcommand{\zp}{\mbox{$\mathbb{Z}_p$}}
\newcommand{\zpos}{\mbox{$\mathbb{Z}_{+}$}}
\newcommand{\zposn}{\mbox{$\mathbb{Z}_{+}^n$}}

\newcommand{\almosteverywhere}{\mbox{$\forall E$}}
\newcommand{\realnumbersextended}{\mbox{$\overline{\mathbb{R}}$}}
\newcommand{\real}{\mbox{$\mathbb{R}$}}
\newcommand{\reals}{\mbox{$\mathbb{R}$}}

\newcommand{\cdiscclosed}{\mbox{$\mathbb{D}_c$}}
\newcommand{\cdiscopen}{\mbox{$\mathbb{D}_o$}}
\newcommand{\cdiscoutopen}{\mbox{$(\mathbb{D}^c)_o$}}
\newcommand{\cminus}{\mbox{$\mathbb{C}^-$}}
\newcommand{\cplus}{\mbox{$\mathbb{C}^+$}}
\newcommand{\im}{\mbox{\rm Im}}
\newcommand{\re}{\mbox{\rm Re}}

\newcommand{\cilspaces}{\mbox{${\rm CILSpaces}$}} 
\newcommand{\cilssp}{\mbox{${\rm CILSP}$}} 
\newcommand{\cilsp}{\mbox{${\rm CILSp}$}} 
\newcommand{\cli}{\mbox{${\rm CLI}$}} 
\newcommand{\climin}{\mbox{${\rm CLI_{min}}$}} 
\newcommand{\cn}{\mbox{$\mathbb{C}^n$}} 
\newcommand{\cmcm}{\mbox{$\mathbb{C}^m$}} 
\newcommand{\linearsubspaces}{\mbox{${\rm LinSubspaces}$}} 
\newcommand{\lat}{\mbox{${\rm Lat}$}} 
\newcommand{\quotientspaces}{\mbox{${\rm Quotientspaces}$}} 
\newcommand{\rinfty}{\mbox{$\mathbb{R}^{\infty}$}}
\newcommand{\rk}{\mbox{$\mathbb{R}^k$}}
\newcommand{\rmrm}{\mbox{$\mathbb{R}^m$}}
\newcommand{\rn}{\mbox{$\mathbb{R}^n$}} 
\newcommand{\rp}{\mbox{$\mathbb{R}^p$}}
\newcommand{\rpos}{\mbox{$\mathbb{R}_{+}$}}
\newcommand{\rposk}{\mbox{$\mathbb{R}_{+}^k$}}
\newcommand{\rposm}{\mbox{$\mathbb{R}_{+}^m$}}
\newcommand{\rposn}{\mbox{$\mathbb{R}_{+}^n$}}
\newcommand{\rposp}{\mbox{$\mathbb{R}_{+}^p$}}
\newcommand{\rr}{\mbox{$\mathbb{R}^r$}}
\newcommand{\rspos}{\mbox{$\mathbb{R}_{s+}$}}
\newcommand{\rsposn}{\mbox{$\mathbb{R}_{s+}^n$}}
\newcommand{\simplexn}{\mbox{$\mathbb{S}_{+}^n$}}
\newcommand{\sn}{\mbox{$\mathbb{S}_{+}^n$}}

\newcommand{\cimod}{\mbox{${\rm CIModules}$}} 

\newcommand{\switch}{\mbox{${\rm Switch}$}} 
\newcommand{\threshold}{\mbox{${\rm Threshold}$}}

\newcommand{\ara}{\mbox{${\rm ARA}$}} 
\newcommand{\dom}{\mbox{\rm Dom}}
\newcommand{\Fc}{\mbox{${\bf F}$}} 
\newcommand{\hara}{\mbox{${\rm HARA}$}} 
\newcommand{\range}{\mbox{\rm Range}}
\newcommand{\rra}{\mbox{${\rm RRA}$}} 
\newcommand{\U}{\mbox{${\rm U}$}} 
\newcommand{\A}{\mbox{${\rm A}$}} 

\newcommand{\fanal}{\mathcal{C}^{\omega}}
\newcommand{\fsmooth}{\mathcal{C}^{\infty}}
\newcommand{\fcont}{\mathcal{C}}
\newcommand{\fpcont}{\mathcal{PC}}
\newcommand{\fpconst}{\mathcal{PC}onst}

\newcommand{\arrow}{\mbox{${\rm arrow}$}}
\newcommand{\blockdiagonal}{\mbox{${\rm Bdiag}$}}
\newcommand{\cnm}{\mbox{$\mathbb{C}^{n \times m}$}}
\newcommand{\cnn}{\mbox{$\mathbb{C}^{n \times n}$}}
\newcommand{\cnp}{\mbox{$\mathbb{C}^{n \times p}$}}
\newcommand{\cpn}{\mbox{$\mathbb{C}^{p \times n}$}}
\newcommand{\cpp}{\mbox{$\mathbb{C}^{p \times p}$}}
\newcommand{\cum}{\mbox{${\rm cum}$}}
\newcommand{\dnn}{\mbox{$\mathbb{R}_{\rm diag}^{n \times n}$}}
\newcommand{\dnnpos}{\mbox{$\mathbb{R}_{\rm +, diag}^{n \times n}$}}
\newcommand{\dnnspos}{\mbox{$\mathbb{R}_{\rm diag, +}^{n \times n}$}}
\newcommand{\imp}{\mbox{${\rm Imprim}$}}
\newcommand{\invertible}{\mbox{${\rm inv}$}}
\newcommand{\rkk}{\mbox{$\mathbb{R}^{k \times k}$}}
\newcommand{\rkm}{\mbox{$\mathbb{R}^{k \times m}$}}
\newcommand{\rkn}{\mbox{$\mathbb{R}^{k \times n}$}}
\newcommand{\rmk}{\mbox{$\mathbb{R}^{m \times k}$}}
\newcommand{\rmm}{\mbox{$\mathbb{R}^{m \times m}$}}
\newcommand{\rmn}{\mbox{$\mathbb{R}^{m \times n}$}}
\newcommand{\rmnarrow}{\mbox{$\mathbb{R}_{\arrow}^{m \times n}$}}
\newcommand{\rmnblockdiagonal}{\mbox{$\mathbb{R}_{{\rm Bdiag}}^{m \times n}$}}
\newcommand{\rmncoordinated}{\mbox{$\mathbb{R}_{{\rm c}}^{m \times n}$}}
\newcommand{\rmr}{\mbox{$\mathbb{R}^{m \times r}$}}
\newcommand{\rnk}{\mbox{$\mathbb{R}^{n \times k}$}}
\newcommand{\rnm}{\mbox{$\mathbb{R}^{n \times m}$}}
\newcommand{\rnn}{\mbox{$\mathbb{R}^{n \times n}$}}
\newcommand{\rnnarrow}{\mbox{$\mathbb{R}_{\arrow}^{n \times n}$}}
\newcommand{\rnncoordinated}{\mbox{$\mathbb{R}_{{\rm c}}^{n \times n}$}}
\newcommand{\rnnspd}{\mbox{$\mathbb{R}_{spd}^{n \times n}$}}
\newcommand{\rnnsspd}{\mbox{$\mathbb{R}_{sspd}^{n \times n}$}}
\newcommand{\rnp}{\mbox{$\mathbb{R}^{n \times p}$}}
\newcommand{\rnr}{\mbox{$\mathbb{R}^{n \times r}$}}
\newcommand{\rpk}{\mbox{$\mathbb{R}^{p \times k}$}}
\newcommand{\rpm}{\mbox{$\mathbb{R}^{p \times m}$}}
\newcommand{\rpn}{\mbox{$\mathbb{R}^{p \times n}$}}
\newcommand{\rpp}{\mbox{$\mathbb{R}^{p \times p}$}}
\newcommand{\rpr}{\mbox{$\mathbb{R}^{p \times r}$}}
\newcommand{\rrn}{\mbox{$\mathbb{R}^{r \times n}$}}
\newcommand{\rrr}{\mbox{$\mathbb{R}^{r \times r}$}}
\newcommand{\sposdefnn}{\mbox{$\mathbb{RSPD}^{n \times n}$}}
\newcommand{\spd}{\mbox{$\mathbb{R}_{\rm spd}^{n \times n}$}}

\newcommand{\dposnn}{\mbox{$\mathbb{R}_{{\rm diag},+}^{n \times n}$}}
\newcommand{\dpossnn}{\mbox{$\mathbb{R}_{{\rm diag},+s}^{n \times n}$}}
\newcommand{\dsposnn}{\mbox{$\mathbb{R}_{{\rm diag},s+}^{n \times n}$}}
\newcommand{\dscnn}{\mbox{$\mathbb{DSC}_{+}^{n \times n}$}}
\newcommand{\dsnn}{\mbox{$\mathbb{DS}_{+}^{n \times n}$}}
\newcommand{\indeximprim}{\mbox{\rm imprim}}
\newcommand{\indexprim}{\mbox{\rm indexprim}}
\newcommand{\imprim}{\mbox{\rm imprim}}
\newcommand{\matrixrowtrunc}{\mbox{$\rm matrowtrunc$}}
\newcommand{\matrixdiagtrunc}{\mbox{$\rm mdtrunc$}}
\newcommand{\mnn}{\mbox{$\mathbb{M}_{+}^{n \times n}$}}
\newcommand{\mposkk}{\mbox{$M_{+}^{k \times k}$}}
\newcommand{\mposmm}{\mbox{$M_{+}^{m \times m}$}}
\newcommand{\mposnn}{\mbox{$M_{+}^{n \times n}$}}
\newcommand{\nn}{\mbox{$\mathbb{N}_n$}}
\newcommand{\permnn}{\mbox{$\mathbb{P}^{n \times n}$}}
\newcommand{\permpp}{\mbox{$\mathbb{P}^{p \times p}$}}
\newcommand{\posR}[2]{\mbox{$R_{+}^{#1\times#2}$}}
\newcommand{\rposik}{\mbox{$\mathbb{R}_{+}^{\infty\times k}$}}
\newcommand{\rposim}{\mbox{$\mathbb{R}_{+}^{\infty\times m}$}}
\newcommand{\rposin}{\mbox{$\mathbb{R}_{+}^{\infty\times n}$}}
\newcommand{\rposiq}{\mbox{$\mathbb{R}_{+}^{\infty\times q}$}}
\newcommand{\rposkk}{\mbox{$\mathbb{R}_{+}^{k \times k}$}}
\newcommand{\rposkm}{\mbox{$\mathbb{R}_{+}^{k \times m}$}}
\newcommand{\rposkn}{\mbox{$\mathbb{R}_{+}^{k \times n}$}}
\newcommand{\rposkp}{\mbox{$\mathbb{R}_{+}^{k \times p}$}}
\newcommand{\rposmm}{\mbox{$\mathbb{R}_{+}^{m \times m}$}}
\newcommand{\rposmk}{\mbox{$\mathbb{R}_{+}^{m \times k}$}}
\newcommand{\rposmn}{\mbox{$\mathbb{R}_{+}^{m \times n}$}}
\newcommand{\rposnk}{\mbox{$\mathbb{R}_{+}^{n \times k}$}}
\newcommand{\rposnm}{\mbox{$\mathbb{R}_{+}^{n \times m}$}}
\newcommand{\rposnn}{\mbox{$\mathbb{R}_{+}^{n \times n}$}}
\newcommand{\rposnp}{\mbox{$\mathbb{R}_{+}^{n \times p}$}}
\newcommand{\rpospm}{\mbox{$\mathbb{R}_{+}^{p \times m}$}}
\newcommand{\rpospn}{\mbox{$\mathbb{R}_{+}^{p \times n}$}}
\newcommand{\rposqm}{\mbox{$\mathbb{R}_{+}^{q \times m}$}}
\newcommand{\rposqn}{\mbox{$\mathbb{R}_{+}^{q \times n}$}}
\newcommand{\rposqq}{\mbox{$\mathbb{R}_{+}^{q \times q}$}}
\newcommand{\rsposnn}{\mbox{$\mathbb{R}_{s+}^{n \times n}$}}
\newcommand{\Table}{\mbox{$\rm table$}}
\newcommand{\wnn}{\mbox{$\W^{n \times n}$}}

\newcommand{\adjoint}{\mbox{$\rm Adj$}}
\newcommand{\affine}{\mbox{$\rm Affine$}}
\newcommand{\affineb}{\mbox{$\rm AffB$}}
\newcommand{\blockdiag}{\mbox{{\rm Block-diag}}}
\newcommand{\circulant}{\mbox{$\rm Circulant$}}
\newcommand{\cols}{\mbox{$\rm col$}}
\newcommand{\diag}{\mbox{$\rm Diag$}}
\newcommand{\dsc}{doubly stochastic circulant }
\newcommand{\dscs}{doubly stochastic circulants }
\newcommand{\eig}{\mbox{$\rm eig$}}
\newcommand{\lspan}{\mbox{\rm span}}
\newcommand{\posr}[1]{#1\mbox{\rm -pos-rank}}
\newcommand{\posrank}{\mbox{$\rm pos-rank$}}
\newcommand{\projection}{\mbox{$\rm proj$}}
\newcommand{\rank}{\mbox{$\rm rank$}}
\newcommand{\rowrankz}{\mbox{$\rm row-rank_Z$}}
\newcommand{\specrad}{\mbox{$\rm specrad$}}
\newcommand{\spec}{\mbox{$\rm spec$}}
\newcommand{\spectrum}{\mbox{$\rm spec$}}
\newcommand{\tr}{\mbox{$\rm tr$}}
\newcommand{\trace}{\mbox{\rm trace}}

\newcommand{\affinehull}{\mbox{$\rm affh$}}
\newcommand{\convexhull}{\mbox{$\rm convh$}}
\newcommand{\ri}{\mbox{\rm ri}}

\newcommand{\cekk}{\mbox{$CE_{k,k}$}}
\newcommand{\cekm}{\mbox{$CE_{k,m}$}}
\newcommand{\cekn}{\mbox{$CE_{k,n}$}}
\newcommand{\cenn}{\mbox{$CE_{n,n}$}}
\newcommand{\ckk}{\mbox{$\mathbb{C}_{k,k}$}}
\newcommand{\ckm}{\mbox{$\mathbb{C}_{k,m}$}}
\newcommand{\cone}{\mbox{\rm cone}} 
\newcommand{\faces}{\mbox{\rm Faces}}
\newcommand{\hp}{\mbox{\rm HyperPlane}}
\newcommand{\normalvector}{\mbox{$V_{normal}$}}
\newcommand{\phs}{\mbox{$PHS$}}
\newcommand{\plsets}{\mbox{$PLSets$}}
\newcommand{\polyhedralcone}{\mbox{\rm Polyhcone}}
\newcommand{\shp}{\mbox{\rm SupportHyperPlane}}
\newcommand{\subpolytope}{\mbox{\rm SubPolytope}}
\newcommand{\subrectangle}{\mbox{\rm SubRectangle}}
\newcommand{\vertices}{\mbox{$V_{vertices}$}}

\newcommand{\qy}{\mbox{$Q_y$}}
\newcommand{\qlsdp}{\mbox{${\bf \partial Q_{lsdp}}$}}
\newcommand{\qp}{\mbox{${\bf Q_{lsp}}$}}
\newcommand{\qpd}{\mbox{${\bf Q_{lsdp}}$}}
\newcommand{\qpr}{\mbox{${\bf Q_{lsp,r}}$}}
\newcommand{\qprs}{\mbox{${\bf Q_{lsp,s}}$}}
\newcommand{\qpdr}{\mbox{${\bf Q_{lsdp,r}}$}}
\newcommand{\dqps}{\mbox{${\bf \partial Q_{lsp,r}}$}}
\newcommand{\dqpss}{\mbox{${\bf \partial Q_{lsp,s}}$}}
\newcommand{\dqpdr}{\mbox{${\bf \partial Q_{lsdp,r}}$}}
\newcommand{\dqprs}{\mbox{${\bf \partial Q_{lsp,r,s}}$}}
\newcommand{\dqpdrs}{\mbox{${\bf \partial Q_{lsdp,r,s}}$}}
\newcommand{\dqp}{\mbox{${\bf \partial Q_{lsp}}$}}
\newcommand{\dqpd}{\mbox{${\bf \partial Q_{lsdp}}$}}

\newcommand{\con}{\mbox{${\rm con}$}}
\newcommand{\conmat}{\mbox{${\rm conmat}$}}
\newcommand{\conset}{\mbox{${\rm Conset}$}}
\newcommand{\coconset}{\mbox{${\rm co-Conset}$}}
\newcommand{\controllablepair}{\mbox{${\rm conpair}$}}
\newcommand{\controllabilitymatrix}{\mbox{${\rm conmat}$}}
\newcommand{\controllableset}{\mbox{${\rm conset}$}}
\newcommand{\cocontrollableset}{\mbox{${\rm co-conset}$}}
\newcommand{\controlset}{\mbox{$\rm controlset$}}
\newcommand{\ls}{\mbox{${\rm LS}$}}
\newcommand{\lsp}{\mbox{${\rm LSP}$}}
\newcommand{\lspmin}{\mbox{${\rm LSP_{min}}$}}
\newcommand{\reachm}{\mbox{${\rm reachm}$}}
\newcommand{\obsm}{\mbox{${\rm obsm}$}}
\newcommand{\obsmat}{\mbox{${\rm obsmat}$}}
\newcommand{\obsmap}{\mbox{${\rm obsmap}$}}
\newcommand{\realization}{\mbox{${\rm realiz}$}}
\newcommand{\reconmap}{\mbox{${\rm reconmap}$}}

\newcommand{\argmin}{\mbox{${\rm argmin}$}}
\newcommand{\argmax}{\mbox{${\rm argmax}$}}
\newcommand{\rai}{\mbox{${\rm RAI}$}}

\newcommand{\obs}{\mbox{${\rm obs}$}}
\newcommand{\aobs}{\mbox{${\rm A_{obs}}$}}
\newcommand{\qobs}{\mbox{${\rm Q_{obs}}$}}
\newcommand{\svdtruncation}{\mbox{${\rm SVDtrunc}$}}
\newcommand{\zclosure}{\mbox{${\rm \mathcal{Z}-cl}$}}

\newcommand{\act}{\mbox{$\rm act$}}
\newcommand{\aux}{\mbox{$\rm Aux$}}
\newcommand{\cat}{\mbox{$\rm cat$}}
\newcommand{\cig}{\mbox{$\rm CIG$}}
\newcommand{\cil}{\mbox{$\rm CIL$}}
\newcommand{\child}{\mbox{$\rm Chi$}}
\newcommand{\co}{\mbox{$\rm CO$}}
\newcommand{\codes}{\mbox{$\rm CODES$}}
\newcommand{\coreach}{\mbox{$\rm coreach$}}
\newcommand{\coreachcomponent}{\mbox{$\rm coreachco$}}
\newcommand{\coreachgen}{\mbox{$\rm coreachgen$}}
\newcommand{\coreachset}{\mbox{$\rm coreachset$}}
\newcommand{\csublanguage}{\mbox{$\rm C$}}
\newcommand{\cpcsublanguage}{\mbox{$\rm C_{pc}$}}
\newcommand{\csuplanguage}{\mbox{$\rm CSupL$}}
\newcommand{\cpcsuplanguage}{\mbox{$\rm CSupL_{pc}$}}
\newcommand{\decdes}{\mbox{$\rm decDES$}}
\newcommand{\infcsuplanguage}{\mbox{$\inf {\rm CSupL}$}}
\newcommand{\infcsuppclanguage}{\mbox{$\inf {\rm CSupL_{pc}}$}}
\newcommand{\last}{\mbox{$\rm last$}}
\newcommand{\length}{\mbox{$\rm length$}}
\newcommand{\markact}{\mbox{$\rm markact$}}
\newcommand{\moddes}{\mbox{$\rm modDES$}}
\newcommand{\nextact}{\mbox{$\rm nextact$}}
\newcommand{\parent}{\mbox{$\rm Par$}}
\newcommand{\prefix}{\mbox{$\rm prefix$}}
\newcommand{\ps}{\mbox{$\rm PS$}}
\newcommand{\qeq}{\mbox{$\rm QEQ$}}
\newcommand{\reachcomponent}{\mbox{$\rm reachco$}}
\newcommand{\reachgen}{\mbox{$\rm reachgen$}}
\newcommand{\reachset}{\mbox{$\rm reachset$}}
\newcommand{\runs}{\mbox{$\rm Runs$}}
\newcommand{\selfloop}{\mbox{$\rm selfloop$}}
\newcommand{\setldfg}{\mbox{$\rm SetL_{DFG}$}}
\newcommand{\setllc}{\mbox{$\rm SetL_{lc}$}}
\newcommand{\setllcn}{\mbox{$\rm SetL_{lcn}$}}
\newcommand{\setln}{\mbox{$\rm SetL_N$}}
\newcommand{\setlnfg}{\mbox{$\rm SetL_{NFG}$}}
\newcommand{\setlreg}{\mbox{$\rm SetL_{reg}$}}
\newcommand{\suffix}{\mbox{$\rm suffix$}}
\newcommand{\shuffle}{\mbox{$\rm shuffle$}}
\newcommand{\size}{\mbox{$\rm Size$}}
\newcommand{\starrr}{\mbox{$\rm Star$}}
\newcommand{\supap}{\mbox{$\sup {\rm AP}$}}
\newcommand{\supc}{\mbox{$\sup {\rm C}$}}
\newcommand{\supcc}{\mbox{$\sup {\rm cC}$}}
\newcommand{\supccn}{\mbox{$\sup {\rm cCN}$}}
\newcommand{\supcsublanguage}{\mbox{$\sup {\rm C}$}}
\newcommand{\supcpcsublanguage}{\mbox{$\sup {\rm C_{pc}}$}}
\newcommand{\supcnsublanguage}{\mbox{$\sup {\rm CN}$}}
\newcommand{\supcn}{\mbox{$\sup {\rm CN}$}}
\newcommand{\supmccn}{\mbox{$\sup {\rm mcCN}$}}
\newcommand{\suppc}{\mbox{$\sup {\rm PC}$}}
\newcommand{\suppn}{\mbox{$\sup {\rm PN}$}}
\newcommand{\supn}{\mbox{$\sup {\rm N}$}}
\newcommand{\TIME}{\mbox{$\rm TIME$}}
\newcommand{\transrel}{\mbox{$\rm Tr$}}
\newcommand{\trim}{\mbox{$\rm trim$}}
\newcommand{\trimgen}{\mbox{$\rm trimgen$}}
\newcommand{\triple}{\mbox{$\rm tri$}}
\newcommand{\uc}{\mbox{$\rm uc$}}

\newcommand{\cogstocsp}{\mbox{${\rm COGStocSP}$}}
\newcommand{\fstocs}{\mbox{${\rm FStocS}$}}
\newcommand{\fss}{\mbox{${\rm FSS}$}}
\newcommand{\fstocsp}{\mbox{${\rm FStocSP}$}}
\newcommand{\gstocs}{\mbox{${\rm GStocS}$}}
\newcommand{\gstocsp}{\mbox{${\rm GStocSP}$}}
\newcommand{\gstoccsp}{\mbox{${\rm GStocCSP}$}}
\newcommand{\stoccs}{\mbox{${\rm StocCS}$}}
\newcommand{\stocs}{\mbox{${\rm StocS}$}}

\newcommand{\is}{\mbox{${\rm IS}$}}
\newcommand{\isdsrs}{\mbox{${\rm ISDSrs}$}}
\newcommand{\isdsones}{\mbox{${\rm ISDS1s}$}}
\newcommand{\isdsrscommon}{\mbox{${\rm ISDSrsCommon}$}}
\newcommand{\isdsonescommon}{\mbox{${\rm ISDS1sCommon}$}}
\newcommand{\isdsrsprivate}{\mbox{${\rm ISDSrsPrivate}$}}
\newcommand{\isnn}{\mbox{${\rm ISNN}$}}
\newcommand{\inn}{\mbox{${\rm I_{nn}}$}}

\newcommand{\bits}{\mbox{$\rm bits$}}

\newcommand{\dtime}{\mbox{$\rm DTIME$}}
\newcommand{\exptime}{\mbox{$\rm EXPTIME$}}
\newcommand{\np}{\mbox{$\rm NP$}}
\newcommand{\ntime}{\mbox{$\rm NTIME$}}
\newcommand{\polylog}{\mbox{$\rm polylog$}}
\newcommand{\timecomplexity}{\mbox{$\rm TIME$}}

\newcommand{\glnr}{\mbox{$Gl_{n}(\mathbb{R})$}}

\newcommand{\aut}{\mbox{${\bf Aut}$}}
\newcommand{\catset}{\mbox{${\bf Set}$}}
\newcommand{\comp}{\mbox{$comp$}}
\newcommand{\Grp}{\mbox{${\bf Grp}$}}
\newcommand{\kaut}{\mbox{${\bf K-Aut}$}}
\newcommand{\kmedv}{\mbox{${\bf K-Medv}$}}
\newcommand{\kmedvio}{\mbox{$({\bf K-Medv} \downarrow <I^+,O>)$}}
\newcommand{\moduler}{\mbox{${\bf Module_R}$}}
\newcommand{\ob}{\mbox{$ob$}}
\newcommand{\Set}{\mbox{${\bf Set}$}}
\newcommand{\cattopo}{\mbox{${\bf Topo}$}}

\newcommand{\cont}{\mbox{${\rm cont}$}}
\newcommand{\Deg}{\mbox{${\rm Deg}$}}
\newcommand{\degree}{\mbox{${\rm deg}$}}
\newcommand{\diff}{\mbox{${\rm diff}$}}
\newcommand{\gmon}{\mbox{${\rm G_{mon}}$}}
\newcommand{\diffrpos}{\mbox{${\rm Diff} \mathbb{R}_+$}}
\newcommand{\lcm}{\mbox{${\rm lcm}$}}
\newcommand{\mon}{\mbox{${\rm mon}$}}
\newcommand{\mnm}{\mbox{${\rm mnm}$}}
\newcommand{\order}{\mbox{${\rm order}$}}
\newcommand{\rpoly}[2]{\mbox{$R_+[#1]/(#1^{#2}-1)$}}
\newcommand{\slsp}{\mbox{${\rm SL}\Sigma{\rm P}$}}
\newcommand{\spoly}[2]{\mbox{$S_+[#1]/(#1^{#2}-1)$}}
\newcommand{\sqfree}{\mbox{${\rm sqfree}$}}
\newcommand{\support}{\mbox{${\rm support}$}}
\newcommand{\trdeg}{\mbox{${\rm trdeg}$}}

\newcommand{\as}{\mbox{{\rm $a.s.$}}} 
\newcommand{\aslim}{\mbox{{\rm $a.s.-\lim$}}} 
\newcommand{\ci}{\mbox{{\rm $CI$}}} 
\newcommand{\ciffg}{\mbox{$(F_1,F_2 | G ) \in \ci}} 
\newcommand{\dlim}{\mbox{{\rm $D-\lim$}}} 
\newcommand{\essinf}{\mbox{{\rm $essinf$}}} 
\newcommand{\esssup}{\mbox{{\rm $esssup$}}} 
\newcommand{\foralmostall}{\mbox{{\rm $\mbox{almost all}$}}} 
\newcommand{\ift}{\mbox{{\rm $ift$}}} 
\newcommand{\ltwolim}{\mbox{{\rm $L_2-\lim$}}} 
\newcommand{\pdf}{\mbox{{\rm $pdf$}}} 
\newcommand{\pessinf}{\mbox{{\rm $P-essinf$}}} 
\newcommand{\pesssup}{\mbox{{\rm $P-esssup$}}} 
\newcommand{\plim}{\mbox{{\rm $P-\lim$}}} 

\newcommand{\aloc}{\mbox{{\rm ${\bf A_{loc}}$}}}
\newcommand{\alocplus}{\mbox{{\rm ${\bf A_{loc}^+}$}}}
\newcommand{\aone}{\mbox{{\rm ${\bf A_1}$}}}
\newcommand{\aplus}{\mbox{{\rm ${\bf A^+}$}}}
\newcommand{\bvar}{\mbox{{\rm $BV$}}}
\newcommand{\bvarc}{\mbox{{\rm $BV^c$}}}
\newcommand{\cadlag}{\mbox{{\rm c\`{a}dl\`{a}g}}}
\newcommand{\DL}{\mbox{{\rm $DL$}}}
\newcommand{\mone}{\mbox{{\rm $M_1$}}}
\newcommand{\monec}{\mbox{{\rm $M_1^c$}}}
\newcommand{\monepos}{\mbox{{\rm $M_{+,1}$}}}
\newcommand{\moneu}{\mbox{{\rm $M_{1u}$}}}
\newcommand{\moneuc}{\mbox{{\rm $M_{1u}^c$}}}
\newcommand{\moneuloc}{\mbox{{\rm $M_{1uloc}$}}}
\newcommand{\moneulocc}{\mbox{{\rm $M_{1uloc}^c$}}}
\newcommand{\mtwo}{\mbox{{\rm $M_2$}}}
\newcommand{\mtwos}{\mbox{{\rm $M_{2s}$}}}
\newcommand{\mtwosc}{\mbox{{\rm $M_{2s}^c$}}}
\newcommand{\mtwosd}{\mbox{{\rm $M_{2s}^d$}}}
\newcommand{\mtwosloc}{\mbox{{\rm $M_{2sloc}$}}}
\newcommand{\mtwoslocc}{\mbox{{\rm $M_{2sloc}^c$}}}
\newcommand{\mtwoslocd}{\mbox{{\rm $M_{2sloc}^d$}}}
\newcommand{\mtwoc}{\mbox{{\rm $M_2^c$}}}
\newcommand{\var}{\mbox{{\rm $Var$}}}
\newcommand{\semim}{\mbox{{\rm $SemM$}}}
\newcommand{\semimc}{\mbox{{\rm $SemM^c$}}}
\newcommand{\semimone}{\mbox{{\rm $SemM_1$}}}
\newcommand{\semimonec}{\mbox{{\rm $SemM_1^c$}}}
\newcommand{\semimspecial}{\mbox{{\rm $SemM_s$}}}
\newcommand{\semimtwo}{\mbox{{\rm $SemM_2$}}}
\newcommand{\semimtwoc}{\mbox{{\rm $SemM_2^c$}}}
\newcommand{\semimloc}{\mbox{{\rm $SemM_{loc}$}}}
\newcommand{\semimlocc}{\mbox{{\rm $SemM_{loc}^c$}}}
\newcommand{\stoppingtimes}{\mbox{{\rm $T_{st}$}}}
\newcommand{\stoppingtimesinfty}{\mbox{{\rm $T_{st \uparrow \infty}$}}}
\newcommand{\submone}{\mbox{{\rm $SubM_1$}}}
\newcommand{\submonec}{\mbox{{\rm $SubM_1^c$}}}
\newcommand{\submonepos}{\mbox{{\rm $SubM_{+,1}$}}}
\newcommand{\submpos}{\mbox{{\rm $SubM_+$}}}
\newcommand{\submposc}{\mbox{{\rm $SubM_+^c$}}}
\newcommand{\supmone}{\mbox{{\rm $SupM_1$}}}

\newcommand{\deficiency}{\mbox{{\rm dfc}}}
\newcommand{\reactionnet}{\mbox{{\rm rnet}}}
\newcommand{\rnet}{\mbox{{\rm $rnet$}}}
\newcommand{\netc}{\mbox{{\rm $net_c$}}}
\newcommand{\nets}{\mbox{{\rm $net_s$}}}

\newcommand{\lane}{\mbox{$\rm lane$}}
\newcommand{\roadnet}{\mbox{$\rm RoadNet$}}
\newcommand{\roadsection}{\mbox{$\rm Secs$}}
\newcommand{\subnet}{\mbox{$\rm SubNet$}}
\newcommand{\subnetin}{\mbox{$\rm SubNet_{in}$}}
\newcommand{\subnetout}{\mbox{$\rm SubNet_{out}$}}
\newcommand{\subnetlink}{\mbox{$\rm R_{link}$}}

\newcommand{\dl}{D_{\lambda}}
\newcommand{\occ}{\mbox{\rm Occ}}

\newcommand{\ta}{\widetilde{\alpha}}
\newcommand{\tb}{\widetilde{\beta}}
\newcommand{\tg}{\widetilde{\gamma}}
\newcommand{\td}{\widetilde{\delta}}
\newcommand{\te}{\widetilde{\varepsilon}}
\newcommand{\tx}{\widetilde{\xi}}
\newcommand{\tp}{\widetilde{p}}
\newcommand{\tm}{\widetilde{M}}
\newcommand{\wt}[1]{\widetilde{#1}}
\newcommand{\wh}[1]{\widehat{#1}}

\newcommand{\TT}{T\!\!\!\! I}
\newcommand{\DD}{D\!\!\!\! I}
\newcommand{\RR}{R\!\!\!\! I}
\newcommand{\CC}{C\!\!\!\! I}
\newcommand{\NN}{N\!\!\!\!\!\! I\;}      
 


\maketitle
\pagestyle{plain}

\begin{abstract}
Rational observers are to be constructed for rational systems
while polynomial observers are to be constructed for polynomial systems.
An observer synthesis procedure is formulated. 
First an output-based rational realization is synthesized 
for the considered rational system. 
Then a perturbation technique creates an observer.
Finite algebraic observability of the rational system
impies the existence of a output-based rational realization.
Several examples of rational observers are provided 
including a polynomial system of which 
the state-space dimension of the polynomial observer is 
strictly higher than that of the corresponding system.
\end{abstract}
\section{INTRODUCTION}\label{sec:intro}
The aim of this paper is to show how, 
for a rational system,
a rational observer can be synthesized.
The synthesis procedure is illustrated by several examples.
\par
In control theory there is a need for observers of systems.
Observers are used to produce estimates or predictions
of values of the output of a system.
\par
There is a large body of literature on observers of nonlinear systems.
The reader is referred to the next section for a brief literature review.
\par
The focus of this paper is on rational systems
which arise in biochemical reaction systems,
in physiological systems, and
in engineering.
These systems have been investigated by the authors in various papers,
\cite{nemcova:schuppen:2009,nemcova:schuppen:2010,nemcova:petreczky:schuppen:2013:siamjco}.
The scope of the investigation is restricted
by imposing the condition that 
for a rational system one wants a rational observer
while for a polynomial system one wants a polynomial observer.
This is a self-imposed restriction but it makes sense
considering the algebraic framework of rational systems.
\par
A procedure for observer synthesis is proposed.
The procedure consists of several steps 
of which the first one is the construction
of an output-based realization
followed later on by an output injection step.
That observers are basically output-based realizations
of the system generating the output is
due to R.E. Kalman, \cite{kalman:1963:proc},
and to the subsequent research on stochastic realization.
\par
In general a rational or a polynomial observer 
may have a higher state-space dimension 
than the dimension of the system.
Example \ref{ex:observerhigherorder} of a polynomial system
provides a polynomial observer of
a strictly higher state-space dimension than the system.
\par
The outline of the paper follows.
The next section provides a problem description and motivation.
Section \ref{sec:rationalsystems}
defines rational systems and shows how to check their observability.
Section \ref{sec:observersynthesis}
provides the procedure for observer synthesis.
The theory on which the procedure is based is 
provided in Section \ref{sec:theoryobserversynthesis}.
The performance of the observer is discussed
in Section \ref{sec:performanceissues}.
Section \ref{sec:examples} shows several examples of observers. 
%
\section{PROBLEM FORMULATION}\label{sec:problemformulation}
The motivations for the synthesis of an observer of a system
are primarily:
(1) the interest in estimation of the state of the observer for example 
if the state is the concentration of a chemical species or
if it is a concentration of a physiological model of a human being;
(2) the starting point for prediction of a time series
for example prediction of traffic flow in a road network,
\cite{wang:schuppen:vrancken:2014:simpat},
or prediction of produced photo-voltaic power of solar panels; 
(3) control based on partial observations,
\cite{jiehuang:2004}; and
(4) the use of observers in system identification, 
\cite{nemcova:petreczky:schuppen:2015:sysid}.
\par
There is an extensive literature on observers of control systems.
The foundation is the publication of the Kalman filter,
\cite{kalman:1960},
followed by the papers of D. Luenberger on an observer of
a linear system,
\cite{luenberger:1964,luenberger:1966}.
By now there is an extensive literature on 
observers of several classes of nonlinear systems,
\cite{isidori:2001}.
Recent books exclusively on observers include
\cite{besancon:2007,nijmeijer:fossen:1999}
while major papers include
\cite{gauthier:kupke:1994,gauthier:kupke:2001,chaves:sontag:2001}.
\par
What is the definition of an observer?
Several definitions for an observer have been proposed. 
In this paper an observer will be based 
on an output-based realization of a nonlinear system,
meaning a realization of which
the state is a function of the output and of the output's derivatives.
From such a realization one can directly construct the observer.
The view point of an observer being based on an
output-based realization of the system,
is inspired by the publications of Kalman on the
stochastic realization theory of Kalman filters,
see \cite{kalman:1963:proc}.
\par
A restriction is imposed on the algebraic form of the observers.
For a polynomial system the search is restricted
to a polynomial observer
and for a rational system the search is restricted
to a rational observer.
It should be clear that this is a self-imposed restriction,
there may exists observers in a wider class of
systems for example in the class of Nash systems
defined by the authors in
\cite{nemcova:petreczky:schuppen:2013:siamjco}.
A consequence of this restriction on the algebraic
form of the observer is that the observer may have
a higher dimension than the corresponding system,
see Example \ref{ex:observerhigherorder} below.
\par
The performance of the observer is investigated.
The initial condition of the observer is related
to the output function and its derivatives
which derivatives are often not directly available.
Therefore, there is a convergence issue 
in case the observer is started in an initial state
which is different from that of an output-based realization.
Stability analysis of observers is difficult
and one has to prove that the prediction of the observer
for the output converges to the observed output.
This performance criterion is discussed but not
completely treated in this short conference paper.
\begin{problem}
{\em Observer Synthesis}.
Consider a rational system without inputs, 
\begin{eqnarray*}
     \frac{dx(t)}{dt}
  & = & f(x(t)), ~~ x(0) = x_0 \in \mathbb{R}^n,\\
      y(t)
  & = & h(x(t)).
\end{eqnarray*}
Synthesize an observer of the form,
\begin{eqnarray*}
      \frac{dx_o(t)}{dt}
  & = & f_o(x_o(t),y(t)), ~~ x_o(0) = x_{o,0} \in \mathbb{R}^{n_o}, \\
      y_o(t)
  & = & h_o(x_o(t)), ~~ \mbox{such that,} \\
      0
  & = & \lim_{s \rightarrow \infty} [y(s) - y_o(s)], ~
        \mbox{and,}
\end{eqnarray*}
\begin{itemize}
\item
if the system is a rational system ($f$ and $h$ rational maps),
then the observer is a rational system
($f_o$ and $h_o$ are rational maps); and
\item
if the system is a polynomial system ($f$ and $h$ are polynomial maps),
then the observer is a polynomial system,
($f_o$ and $h_o$ are polynomial maps).
\end{itemize}
\end{problem}
Realization theory of discrete-time polynomial systems
was formulated by E.D. Sontag in his Ph.D. thesis,
\cite{sontag:1979:phdthesis},
and generalized by Z. Bartosiewicz to continuous-time
polynomial systems, \cite{bartosiewicz:1988}.
Realization theory of continuous-time nonlinear and rational systems
was initiated in
\cite{jakubczyk:1980,yuanwang:sontag:1992:siamjco} and further developed by
by J. N\v{e}mcov\'{a},
\cite{nemcova:schuppen:2009,nemcova:schuppen:2010}.
\section{RATIONAL SYSTEMS}\label{sec:rationalsystems}
The concepts of rational systems and their observability recalled in this section
are adopted from \cite{bartosiewicz:1988,nemcova:schuppen:2009,nemcova:schuppen:2010}. 
To state the proper definitions,
let us first provide a short overview of necessary terms 
from commutative algebra and algebraic geometry. 
For more details see e.g. 
\cite{zariski:samuel:1958,cox:little:oshea:1992,kunz:1985,bochnak:coste:roy:1998}.
\par
By $\mathbb{R}[X_1, \dots ,X_n]$ we denote 
the algebra of polynomials in $n$ variables 
with coefficients in the real numbers. 
A subset $X \subseteq \mathbb{R}^n$ is called 
a \emph{variety} if it is a set of points of $\mathbb{R}^n$
which satisfy finitely many polynomial equalities. 
We say it is \emph{irreducible} if it cannot be written as a
union of two disjoint varieties. 
Let $I(X)$ denote the ideal of polynomials of 
$\mathbb{R}[X_1, \dots ,X_n]$ vanishing on $X$. 
Then the elements of $\mathbb{R}[X_1, \dots ,X_n]/I(X)$ 
are referred to as \emph{polynomials} on $X$. 
The ring of all such polynomials is denoted by $A_X$. 
Since $X$ is irreducible and thus $A_X$ is an
integral domain, one can define $Q_X$, 
the field of \emph{rational functions} on $X$, as a field of fractions of $A_X$.
\begin{definition}\label{def:rationalsystem}
By a {\em rational system} $\Sigma$ (without inputs) 
we refer to a control system as understood in system theory,
\cite{sontag:1998:book}, with the notation,
\begin{eqnarray}
       \Sigma 
  & =  & (X, Y, f, h, x_0),  \nonumber \\
       \frac{dx(t)}{dt}  
  & = & f(x(t)), ~ x(0) = x_0 \in X, \\
      y(t)  
  & = & h(x(t)),
\end{eqnarray}
where the state-space $X$ is an irreducible variety in $\mathbb{R}^n$, 
$Y = \mathbb{R}^{m_y}$ with $m_y \in \mathbb{Z}_+$, and 
the components of 
$f: X \rightarrow \mathbb{R}^n$ and $h: X \rightarrow \mathbb{R}^{m_y}$ 
are rational maps on $X$ defined at $x_0$.
\par
Polynomial systems are defined analogically 
with the components of $f$ and $h$ being polynomial maps on $X$.
\end{definition}
\begin{example}\label{ex:rationalsystem}
Rational systems are widely used, among others, 
for mathematical description of biological phenomena.
One such example is the following rational system $\Sigma$ 
which describes an enzyme catalyzed change of a substrate to a product. 
The structure of $\Sigma$ is derived by considering the corresponding single
reversible reaction to be modeled by Michaelis-Menten kinetics. 
Let $x_1$ denote the substrate concentration and let $x_2$
denote the product concentration, 
then $\Sigma$ is described by the equations,
\begin{eqnarray*}
      \frac{dx_1(t)}{dt}  
  & = & - a x_1(t) + \frac{c x_1(t) + b x_1^2(t)}{x_1(t) + d}, ~ x_1(0) = 1, \\
      \frac{dx_2(t)}{dt}  
  & = & \frac{e x_1(t)}{x_1(t) + d}, ~ x_2(0) = 1, \\
      y(t)  
  & = & x_2(t),
\end{eqnarray*}
where the considered state-space $X$ equals $\mathbb{R}^2$, 
the concentration $x_2$ of the product is assumed to be observed, 
the initial conditions are chosen to be positive, and 
constant parameters $a, b, c, d, e$ have specific biological meaning.
Let us, for simplicity, assume $a = b = c = e = 1$ and $d = 2$.
\end{example}
In the remainder of the paper the concept of observability
of a rational system is needed which is thus borrowed from
\cite{nemcova:schuppen:2009,nemcova:schuppen:2010}.
\begin{definition}\label{def:observationalgebra}
Consider a rational system $\Sigma = (X, Y, f, h, x_0)$. 
Define the {\em observation algebra} $A_{obs}(\Sigma)$ of $\Sigma$
as the algebra over the real numbers generated by 
the components of $h$ and 
closed with respect to Lie derivatives along the vector field $f$. 
Thus,
\begin{equation}
  A_{obs}(\Sigma)~ 
  = ~ \mathbb{R}[ \{L_f^k h_i ~|~i = 1, \dots,m, ~k =0, 1, \dots \}],
\end{equation}
where $L_f^0 h_i = h_i$, 
$ L_f^1 h_i = L_f h_i= \sum_{j=1}^n f_j(x) \frac{\partial}{\partial x_j} h_i$,
and 
$L_f^k h_i = L_f (L_f^{k-1} h_i)$ for $k =  2, 3 \dots$, $i = 1, \dots, m$.
Because $X$ was an irreducible variety, 
one can define the {\em observation field} $Q_{obs}(\Sigma)$ of $\Sigma$
as the field of fractions of $A_{obs}(\Sigma)$, i.e.
\begin{equation}
  Q_{obs}(\Sigma)~ 
  = ~ \{ p/q | ~ p, ~ q \in A_{obs}(\Sigma), ~ q \neq 0 \}.
\end{equation}
\par
The rational system $\Sigma$ is called {\em algebraically observable}
if its observation field equals the field of all rational functions 
on the state-space, $Q_X = Q_{obs}(\Sigma)$.
\par
The polynomial system $\Sigma$ is called {\em algebraically observable}
if its observation algebra equals 
the algebra of all polynomials on the state-space, $A_X = A_{obs}(\Sigma)$.
\end{definition}
\begin{example}\label{ex:observabilityratsys}
Consider the rational system defined in Example \ref{ex:rationalsystem}. 
Let us refer to it as to $\Sigma = (X,Y,f,h,x_0)$.
Then,
\begin{equation*}
   L_f h(x_1,x_2) 
  = \frac{-x_1}{x_1 + 2} \frac{\partial h(x_1,x_2)}{\partial x_1} 
    + \frac{x_1}{x_1 + 2} \frac{\partial h(x_1,x_2)}{\partial x_2}
\end{equation*}
and $h(x_1,x_2)  =  x_2$. 
Let us compute few elements of $Q_{obs}(\Sigma)$:
\begin{eqnarray*}
     h 
  & = & x_2 \in Q_{obs}(\Sigma),\\
      L_f \left( x_2 \right) 
  & = & \frac{x_1}{x_1 + 2} \in Q_{obs}(\Sigma),\\
      L_f \left( \frac{x_1}{x_1 + 2} \right) 
  & = & \frac{- 2 x_1}{(x_1 + 2)^3} \in Q_{obs}(\Sigma).
\end{eqnarray*}
To obtain further elements one can keep on 
computing Lie derivatives of known elements, 
multiply and divide the elements and multiply them by real numbers.
For example, 
$\frac{- 2 x_1}{(x_1 + 2)^3} \frac{-1}{2} =  \frac{x_1}{(x_1 + 2)^3} \in Q_{obs}(\Sigma)$ 
and 
$\frac{x_1}{(x_1 + 2)} \frac{(x_1 + 2)^3}{x_1} = (x_1 +2)^2 \in Q_{obs}(\Sigma)$.
Then, $L_f \left( (x_1 +2)^2 \right) = -2 x_1$ and 
consequently $x_1$ belong to $Q_{obs}(\Sigma)$. 
Because $x_1, x_2 \in Q_{obs}(\Sigma)$, 
it follows that $Q_{obs}(\Sigma) = \mathbb{R}(x_1,x_2) = Q_X$ 
which implies $\Sigma$ is algebraically observable.
\end{example}
\section{OBSERVER SYNTHESIS}\label{sec:observersynthesis}
The observer synthesis procedure is stated below
and is complemented in the next section 
by theory of the various steps.
Below either a rational or a polynomial system is considered.
A polynomial system is taken along in the exposition
to illustrate the various steps.
\begin{procedure}\label{proc:observersynthesis}
{\em Observer Synthesis Procedure}
\begin{enumerate}
\item
Construct a state-space transformation\\
$s: X \rightarrow \hat{X}$.
\item
Prove the existence of an inverse function $s^{-1}$ of $s$
of the state-space transformation
in the subclass of rational or of polynomial functions
depending on the case.
\item
Derive an output-based realization.
\item
Derive the observer.
\item
Choose an observer gain to meet performance objectives, 
including stability. 
\end{enumerate}
\end{procedure}
\begin{example}\label{ex:polsyspolobs1}
{\em Polynomial system}.
Consider the observable polynomial system
for which an observer is wanted.
\begin{eqnarray*}
     \frac{dx(t)}{dt}
  & = & \left(
        \begin{array}{l}
          - a_{11} x_1(t)^3 + a_{12} x_2(t) \\
          - a_{22} x_2(t)
        \end{array}
        \right), ~ x(0) = x_{0}, \\
      y(t)
  & = & x_1(t) = C x(t) =
        \left(
        \begin{array}{ll}
          1 & 0
        \end{array}
        \right) x(t), \\
  &   & a_{11}, ~ a_{22} \in (0,\infty), ~
        a_{12} \in \mathbb{R}, ~
        a_{12} \neq 0.\\
\end{eqnarray*}
\end{example}
\begin{procedure}\label{proc:step1}
{\em Step 1. Construction of a state-space transformation}.
Define,
\begin{eqnarray*}
      s_1(x)
  & = & h(x),  ~
        s_1: X \rightarrow \mathbb{R}^{m_y}, ~~
        y(t) = h(x(t)), \\
      \frac{dy(t)}{dt}
  & = & \frac{dh(x(t))}{dt}
        = \frac{\partial h(x)}{\partial x}|_{x=x(t)} \frac{dx(t)}{dt} \\
  & = & \frac{\partial h(x)}{\partial x} f(x)|_{x=x(t)} 
        = s_2(x(t)), \\
      s_2(x)
  & = & \frac{\partial s_1(x)}{\partial x} f(x), ~~
        s_2: X \rightarrow \mathbb{R}^{m_y}, \\
      s_{k+1}(x)
  & = & \frac{\partial s_k(x)}{\partial x} f(x), ~
        k \in \mathbb{Z}_+, \\
  &   & m_o \in \mathbb{Z}_+, ~
        n_o = m_o \times m_y, \\
      s(x)
  & = & \left(
        \begin{array}{l}
          s_1(x) \\ s_2(x) \\ \vdots \\ s_{m_o}(x)
        \end{array}
        \right), ~
        s: X = \mathbb{R}^n  \rightarrow \mathbb{R}^{n_o}.
\end{eqnarray*}
\end{procedure}
The choice of the rational observability index $m_o \in \mathbb{Z}_+$  
is best made in Step 2 based on 
the existence of the inverse of the state-space transformation $s$ 
in a specified set. 
\begin{example}\label{ex:polsyspolobs2}
{\em Polynomial system}.
Example \ref{ex:polsyspolobs1} is continued.
The state-space transformation function is calculated.
\begin{eqnarray*}
      s_1(x)
  & = & h(x) = x_1, ~~
      s_2(x)
      = - a_{11} x_1^3 + a_{12} x_2, \\
      s(x)
  & = & ( s_1(x), ~ s_2(x))^T.
\end{eqnarray*}
\end{example}
\begin{procedure}\label{proc:step2}
{\em Step 2. Prove the existence and construct 
a rational inverse or a polynomial inverse 
of the state-space transformation},
\begin{eqnarray*}
      \hat{x}(t)
  & = & \left(
        \begin{array}{llll}
          y(t) & \frac{dy(t)}{dt}  & \ldots 
               & \frac{d^{m_o-1}y(t)}{dt^{m_o-1}}
       \end{array}
       \right)^T, \\
  &   & \hat{x}: T \rightarrow \mathbb{R}^{n_o}, \\
      \hat{x}(t)
  & = & s(x(t)), ~
        s \in Q_{obs}(\Sigma)^{n_o}. 
\end{eqnarray*}
Construct an inverse function $s^{-1}$ of the state-space
transformation $s$ such that
if $s$ is a rational, 
then $s^{-1}$ will be a rational 
while if $s$ is a polynomial, 
then $s^{-1}$ will be a polynomial.
\par
It is proven in the next section that 
observability of the rational or polynomial system 
implies the existence of the inverse 
in the required class of algebraic objects.
\end{procedure}
\begin{example}\label{ex:polsyspolobs3}
{\em Polynomial system}.
Example \ref{ex:polsyspolobs2} is continued.
The state-space transformation $s$ admits a
polynomial inverse,
\begin{eqnarray*}
     \hat{x}_1
  & = & s_1(x) = x_1, ~
      x_1
      = \hat{x}_1, \\
      \hat{x}_2
  & = & s_2(x) = - a_{11} x_1^3 + a_{12} x_2, \\
      s^{-1}(\hat{x})
  & = & \left(
        \begin{array}{l}
          x_1 \\ x_2
        \end{array}
        \right) =
        \left(
        \begin{array}{l}
          \hat{x}_1 \\
          \frac{a_{11}}{a_{12}} \hat{x}_1^3 + \frac{1}{a_{12}} \hat{x}_2
        \end{array}
        \right).
\end{eqnarray*}
\end{example}
\begin{procedure}\label{proc:step3}
{\em Step 3. Derivation of the output-based realization}.
Next the output-based realization can be calculated,
\begin{eqnarray*}
      x(t) 
  & = & s^{-1}(\hat{x}(t)), \\
      \frac{d \hat{x}(t)}{dt}
  & = & f_{or}(\hat{x}(t)), \\
      y(t)
  & = & h_{or}(\hat{x}(t)) = C_o \hat{x}(t) = \hat{x}_{1:m_y}(t), \\
      f_{or}(\hat{x})
  & = & \frac{\partial s(x)}{\partial x} f(x)|_{x = s^{-1}(\hat{x})}, \\
      C_o
  & = & \left(
        \begin{array}{llll}
          I_{m_y} & 0 & \ldots & 0 
        \end{array}
        \right).
\end{eqnarray*}
\end{procedure}
\begin{definition}\label{def:singleoutputcase}
{\em Special case of single output}.
Below attention is restricted to a system
with a one-dimensional output, hence $m_y = 1$.
Then $n_o = m_o$.
Note that then,
\begin{eqnarray*}
      \hat{x}(t)
  & = & \left(
        \begin{array}{llll}
          y(t) & \frac{dy(t)}{dt} & \ldots 
               & \frac{d^{n_o-1}y(t)}{dt^{n_o-1}}
       \end{array}
       \right)^T, \\
      \frac{d}{dt} \hat{x}_i(t)
  & = & \frac{d}{dt} \frac{d^{i-1} y(t)}{dt^{i-1}} 
        = \frac{d^{i}y(t)}{dt^{i}}
        = \hat{x}_{i+1}(t), ~~\\
  &   & i = 1, 2, \ldots, n_o-1, \\
      \frac{d}{dt} \hat{x}_{n_o}(t)
  & = & \frac{d}{dt} \frac{d^{n_o-1}y(t)}{dt^{n_o-1}} 
        = \frac{d^{n_o}y(t)}{dt^{n_o}} 
        = s_{n_o+1}(x(t)). \\
      \frac{d\hat{x}(t)}{dt}
  & = & f_{or}(\hat{x}(t)) 
        = A_o \hat{x}(t) + b_o(\hat{x}(t)), ~~\\
      A_o
  & = & \left(
        \begin{array}{lllll}
          0      & 1 & 0      & \ldots  & 0 \\
          0      & 0 & 1      & \ldots  & 0 \\
          \vdots &   & \ddots &         & \vdots \\
          0      & 0 & 0      & \ldots  & 1 \\
          0      & 0 & 0      & 0       & 0
        \end{array}
        \right), ~\\
  &   & b_o(\hat{x}) = 0, \quad \mbox{except that,} \\
      b_{o,n_o}(\hat{x})
  & = & s_{n_o+1}(s^{-1}(\hat{x}))
        = \frac{\partial s_{n_o}(x)}{\partial x} f(x)|_{x = s^{-1}(\hat{x})}.
\end{eqnarray*}
The same structure holds in the multi-output case,
with $m_y > 1$ though then the matrix $A_o$ and the vector $b_o$
have multivariable components.
\end{definition}
%
\begin{example}\label{ex:polsyspolobs4}
{\em Polynomial system}.
Example \ref{ex:polsyspolobs3} is continued.
The output-based realization is,
\begin{eqnarray*}
      \frac{d \hat{x}(t)}{dt}
  & = & \left(
        \begin{array}{l}
          \hat{x}_{2}(t) \\
          \hat{f}_{or,2}(\hat{x}(t))
        \end{array}
        \right), ~
        \hat{x}(0) = s(x(0)), \\
      \hat{f}_{or,2}(\hat{x})
  & = & - a_{11} a_{22} \hat{x}_1^3
        - 3 a_{11} \hat{x}_1^2 \hat{x}_2 
        - a_{22} \hat{x}_2, \\
      y(t)
  & = & C_{or} \hat{x}(t) = \hat{x}_1(t).
\end{eqnarray*}
\end{example}
In an output-based realization, 
the output is a component of the state vector of the realization.
What is needed is that the output of the system becomes
an input of the observer.
This can be achieved by replacing the first component
of the output-based realization by an abstract variable
and by injecting the output into the first component.
Due to the fact that in the output-based realization
the output function is linear,
the injection of the output in the observer is with the
linear function $y(t) - C_o x_o(t)$.
\begin{procedure}\label{proc:step4}
{\em Step 4. Derivation of the observer by output injection}.
Define the observer as the system,
\begin{eqnarray*}
      \frac{dx_o(t)}{dt}
  & = & A_o x_o(t) + b_o (x_o(t)) + \\
  &   & + [ k_o(x_o(t)) + K] [ y(t) - C_o x_o(t)], \\
  &   & x_o(t_0) = x_{o,0} \in \mathbb{R}^n, ~
        K \in \mathbb{R}^{n_0 \times m_y}, ~\\
      k_o(x_o) 
  & = & \left(
        \begin{array}{llll}
          0 & 0 & \ldots & k_{o, m_o}(x_o)
        \end{array}
        \right)^T \in \mathbb{R}^{n_o \times m_y}, \\
      k_{o,m_o}(x_o)
  & = & \frac{\partial b_{o,m_o}(x_o)}{\partial x_{o,1}}
          \in \mathbb{R}^{m_y \times m_y}.
\end{eqnarray*}
\end{procedure}
An explanation of the derivation of the observer follows
for the case $m_y = 1$.
Apply a Taylor expansion to the output-based realization,
and retain only the first-order term, 
\cite{ortega:rheinboldt:1970}.
\begin{eqnarray*}
      \frac{d \hat{x}(t)}{dt}
  & = & f_{or}(\hat{x}(t)), ~~ \hat{x}(0) = \hat{x}_0 \in \mathbb{R}^{n_o},\\
      y(t)
  & = & h_{or}(\hat{x}(t)), \\
      \hat{x}(t)
  & = & \left(
        \begin{array}{lllll}
          \hat{x}_1(t) & \hat{x}_2(t) & \hat{x}_3(t) & \ldots & \hat{x}_{n_o}(t)
        \end{array}
        \right)^T, \\
      \overline{x}(t)
  & = & \left(
        \begin{array}{lllll}
          y(t) & \hat{x}_2(t) & \hat{x}_3(t) & \ldots & \hat{x}_{n_o}(t)
        \end{array}
        \right)^T, \\
      \lefteqn{
        f_{or}(\hat{x}(t))
      } \\
  & = & f_{or}(\overline{x}(t))
        + [ f_{or}(\hat{x}(t)) - f_{or}(\overline{x}(t)) ], \\
  & \approx & f_{or}(\overline{x}(t))
        + \frac{\partial f_{or}(x)}{\partial x}|_{x = \overline{x}(t)}
          [ \hat{x}(t) - \overline{x}(t) ] \\
  & = & f_{or}(\overline{x}(t))
        + \frac{\partial f_{or}(x)}{\partial x_1}|_{x = \overline{x}(t)}
          [ y(t) - C_o \overline{x}(t)) ], \\
      \frac{d x_o(t)}{dt}
  & = & f_o(x_o(t)) 
        + [ k_o(x_o(t)) + K ] [y(t) - C_o x_o(t) ].
\end{eqnarray*}
The choice for a gain matrix $K \in \mathbb{R}^{n_o \times m_y}$
made above is expected to guarantee local stability for specific values.
Further research is required on 
the global asymptotic stability of the performance system, 
see Section VI, to allow the formulation of nonlinear gain functions.
\begin{example}\label{ex:polsyspolobs5}
{\em Polynomial system}.
Example \ref{ex:polsyspolobs4} is continued.
\begin{eqnarray*}
     \frac{d x_o(t)}{dt}
  & = & \left(
        \begin{array}{l}
          x_{o,2}(t) \\
          f_{o,2}(x_o(t),y(t))
        \end{array}
        \right) +  \\
  &   & + [ k_o(x_o(t)) + K] [ y(t) - C x_o (t) ] \\
  & = & \left(
        \begin{array}{l}
          x_{o,2}(t) \\
          - a_{11} a_{22} x_{o,1}^3
          - 3 a_{11} x_{o,1}^2 x_{o,2}
          - a_{22} x_{o,2}
        \end{array}
        \right) \\
  &   & +
        \left(
        \begin{array}{l}
        k_1 \\
        k_2 
        - 3 a_{11} a_{22} x_{o,1}^2
        - 6 a_{11} x_{o,1}  x_{o,2}
        \end{array}
        \right) \times \\
  &   & \times
        [ y(t) - C x_o (t) ], \nonumber  \\
      y_o(t)
  & = & C x_o(t).
\end{eqnarray*}
\end{example}
\begin{procedure}\label{proc:step5}
{\em Step 5. Choose the oberver gain}.
Choose the gain matrix such that
the performance system of Def. \ref{def:performancesystem}
meets the performance objective of stability
and good transient response.
\end{procedure}
\section{THEORY OF OBSERVER SYNTHESIS}\label{sec:theoryobserversynthesis}
This section provides concepts and results which show that
the synthesis procedure of the previous section 
produces indeed a rational or a polynomial observer.
\par
Consider a rational system.
Def. \ref{def:observationalgebra}
defined first the observation algebra of a rational system
as the algebra of rational functions
generated by the infinite sequence of Lie derivatives of the system
for the output components and
subsequently the observation field as the field of fractions
of this observation algebra. 
Below a corresponding concept is defined for
a finite sequence of these Lie derivatives of output components.
Note that the finite set of the zero-th upto the $m$th Lie derivative
of the components of the output function
are precisely the family of $\{s_i, ~ i \in \mathbb{Z}_{m+1}\}$ functions.
\par
The exposition below was first written for the single
output case, thus for $m_y = 1$.
It can be read also for the multi-output case,
for $m_y \geq 2$, if the sequence of $s$ functions
and their components are renumbered.
For the rest of the section, renumber the functions,
$s_{1,1}, s_{1,2}, \dots, s_{1,m_y}, s_{2,1}, \dots, s_{2,m_y}, s_{3,1}, \dots$ 
as 
$s_1, s_2, \dots$, i.e. 
$s_1 = s_{1,1}$, $s_2 = s_{1,2}$, $\dots$, 
$s_{m_y} = s_{1,m_y}$, $s_{m_y + 1} = s_{2,1}$, $\dots$.
\begin{definition}\label{def:finiteobservationalgebra}
{\em Finite algebraic observability of a rational and of a polynomial system}.
\begin{itemize}
\item[(a)]
Consider a rational system $\Sigma = (X, Y, f, h, x_0)$
and let the functions 
$\{ s_i:\mathbb{R}^n \rightarrow \mathbb{R}, i \in \mathbb{Z}_+\}$
be as constructed in Procedure \ref{proc:step1}.
Define the 
{\em $m$-th observation algebra} of rational functions
and the associated 
{\em $m$-th observation field}
for all $m \in \mathbb{Z}_+$ as respectively,
\begin{eqnarray*}
      A_{obs,m}(\Sigma)
  & = & \mathbb{R}[\{ s_1, s_2, \ldots, s_m \}], \\
      Q_{obs,m}(\Sigma)
  & = & \{ p / q | ~ p,q \in A_{obs,m}(\Sigma), ~ q \neq 0 \}.
\end{eqnarray*}
Call the rational system
{\em finitely algebraically observable} if
there exists an integer $m \in \mathbb{Z}_+$ such that
$Q_{obs,m}(\Sigma) = \mathbb{R}(x)$.
\par
Call the least-integer $m \in \mathbb{Z}_+$ such that 
$Q_{obs,m}(\Sigma) = \mathbb{R}(x)$,
the {\em rational observability index} of the rational system,
denote it by $m_o \in \mathbb{Z}_+$, and 
let $n_o = m_o * m_y$ (with disregard of renumbering).
\item[(b)]
Consider a polynomial system.
Define correspondingly for this polynomial system
the {\em $m$-th observation algebra} of polynomial functions,
{\em finite algebraic observability}, and
the {\em polynomial observability index}.
\end{itemize}
\end{definition}
The concept of a rational observability index
and its polynomial analogue
are analogous to those for linear systems, see
\cite[p. 356--357]{kailath:1980:book}.
Less useful for observers because of robustness reasons
seems the concept of rational Kronecker indices,
see \cite{murota:1987} for those of linear systems.
Each such index for a component of the output function
is defined as the least integer for which the 
rational observation field of that component stops increasing forever.
The value of these indices depends on 
the ordering of the output components.
\par
It follows from \cite[Prop. 5.7]{nemcova:schuppen:2009}
that if a rational system is algebraically observable
then its observation field is finitely generated
hence the rational system is finitely algebraically observable.
In the case of a polynomial system it is not known
whether algebraic observability implies that this polynomial
system is finitely algebraically observable.
Therefore the finite algebraic observability will be assumed.
\par
Finite algebraic observability of a rational or of a polynomial system
does not imply 
the existence of an inverse of $s$ with $m \geq \lceil n/m_y \rceil$.
Example \ref{ex:observerhigherorder} of a polynomial system
shows that the polynomial observability index 
can be strictly higher than the state-space dimension of the system, 
or, equivalently, with $m_y = 1$, $n_o > n$.
\begin{proposition}\label{prop:existenceratinverse}
{\em Algebraic characterization 
of the existence of a rational or of a polynomial inverse 
of the state-space transformation.}
Consider a polynomial system $\Sigma = (X, Y, f, h, x_0)$,
the functions $\{s_i, i \in \mathbb{Z}_+\}$,
and the family of observation algebras
$\{ A_{obs,m}(\Sigma), ~ m \in \mathbb{Z}_+\}$.
\begin{eqnarray*}
  &   & \exists ~ m \in \mathbb{Z}_+ ~ 
        \exists ~ r_1, r_2, \ldots, r_n \in \mathbb{R}[S_1, \ldots, S_m], \\
  &   & \mbox{and let} ~ s = (s_1, \ldots, s_m), \\
  &   & \mbox{such that}, ~
        x_i = r_i(s(x)), ~ 
        \forall ~ i \in \mathbb{Z}_n; ~ \\
  &   & \mbox{hence,} ~
        x = r(s(x)) \in A_{obs,m}(\Sigma)^n ~ 
        \mbox{and} ~
        s^{-1} = r; \\
  &   & \mbox{if and only if,}  \\
  &   & \exists ~ m \in \mathbb{Z}_+ ~
        \mbox{such that} ~
        A_{obs,m}(\Sigma) = \mathbb{R}[x].
\end{eqnarray*}
Note that if there exists an $r$ as above, 
then $s$ is invertible with the polynomial inverse $s^{-1} = r$.
\par
A corresponding result holds for the 
existence of a rational inverse $s^{-1}$ of $s$
in terms of the observation field $Q_{obs,m}(\Sigma)$.
\end{proposition}
{\bf Proof.}
$\Leftarrow$ 
Because $\mathbb{R}[x] \subseteq A_{obs,m}(\Sigma)$,
for any $i \in \mathbb{Z}_n$ 
there exists a polynomial 
$r_i(s(x)) \in \mathbb{R}[x_1, \ldots, x_n]$
such that $x_i = r_i(s(x))$.
Define $r = (r_1, \ldots, r_n)$.
Then,
\[
  x = r(s(x))
    =
    \left(
    \begin{array}{l}
      r_1(s(x)) \\ \vdots \\ r_n(s(x)) 
    \end{array}
    \right) \in (A_{obs,m}(\Sigma))^n.
\]
From this follows that $s^{-1} = r \in \mathbb{R}[S_1, \ldots, S_m]$
and thus that the inverse $s^{-1}$ of $s$ exists and 
that it has polynomial components.\\
$\Rightarrow$
If an inverse function $r$ in the indicated set exists 
then $x_i = r_i(s(x)) \in A_{obs,m}(\Sigma)$ for all $i \in \mathbb{Z}_n$.
Thus $\mathbb{R}[x] \subseteq A_{obs,m}(\Sigma) \subseteq \mathbb{R}[x]$
where the last inclusion relation is by definition of $A_{obs,m}(\Sigma)$.
Hence equality holds.
\hfill$\square$
\par\vspace{1\baselineskip}\par\noindent
\begin{theorem}\label{th:existencerationaloutputbasedrealization}
{\em Existence of an output-based realization}.
Consider a rational or a polynomial system.
If the system is finitely algebraically observable 
then there exists an output-based realization 
which is rational or polynomial respectively.
Thus, $\Sigma_{or} = (X_{or}, Y, f_{or}, h_{or}, x_{or,0})$,
as specified in Procedure \ref{proc:step3},
is a rational or a polynomial system respectively.
\end{theorem}
{\bf Proof.}
The case of a rational system is considered,
the case of a polynomial system is similar.
The definition of the state vector $\hat{x}$
as a function of the output $y$ and several of its derivatives
makes clear that the realization is output based.
Note that finite algebraic observability of the system 
implies that there exists 
a finite number $m_o \in \mathbb{Z}_+$ such that
$Q_{obs,m_o}(\Sigma) = \mathbb{R}(x)$.
From Proposition \ref{prop:existenceratinverse}
then follows that the function
$s: \mathbb{R}^n \rightarrow \mathbb{R}^{n_o}$
has an inverse function $s^{-1}$ which is also rational.
From the formulas of Procedure \ref{proc:step3} Step 3,
then follows that the functions $(f_{or}, ~ h_{or})$
are rational.
In more detail,
the inverse function $s^{-1}$ is a rational function,
$f_{or}(\hat{x}) = (\partial s(x)/\partial x) f(x)|_{x = s^{-1}(\hat{x})}$,
and the class of rational functions is closed with respect to substitution,
thus $f_{or}$ is a rational function.
Hence the output-based realization is a rational system.
\hfill$\square$
\par\vspace{1\baselineskip}\par\noindent
The existence of a polynomial inverse of a polynomial map
has been characterized for which purpose we quote
the following results from 
\cite{vandenessen:2000}. 
\par
Let $k$ be an arbitrary field of characteristic zero
(note that for our problem $k = \mathbb{R}$ meets this condition).
A necessary condition for a polynomial map
$F \in k[X_1, \ldots, X_n]$ to be invertible
is that $F$ satisfies the {\em Jacobi condition},
that $\det ( J_F) \in k^*$,
where $k^*$ denotes the units of $k$.
In general, the Jacobi condition is not sufficient.
In case $n=1$, one can prove the sufficiency.
For $n \geq 2$ the problem whether $F \in \mathbb{R}[X_1, \ldots, X_n]$
satisfying the Jacobi condition is invertible is known
as the {\em Jacobi conjecture}.
\par
By adding additional assumptions on polynomial maps
satisfying the Jacobi condition,
one can obtain sufficient conditions of the following form.
\begin{theorem}\label{th:existencepolynomialinverse}
\cite[Th. 2.2.16, p. 53]{vandenessen:2000}.
Let $F \in \mathbb{R}[X_1, \ldots, X_n]$ be such that
$\det (J_F) \in k^*$.
If $k[X]$ is integral over $k[F]$ or if the field extension
of $k(F) \subset k(X)$ is Galois,
then $F$ is invertible.
\end{theorem}
For polynomial maps
$F \in \mathbb{R}[X_1, \ldots, X_n]$ with $k$
being an arbitrary field,
A. van den Essen derived a criterion
based on the theory of Gr\"{o}bner bases
which not only decides whether $F$ is invertible
but also provides a procedure to calculate the inverse 
if it exists.
For the respective theorem (more general than quoted here),
see \cite[ Theorem 3.2.1, p.64]{vandenessen:2000}
or \cite{vandenessen:1990}.
\par
Note that for $n=1$ the invertible polynomials
$F \in k[X]$ are only the affine ones,
for example $F(x) = a x + b$.
For $n=2$ every invertible polynomials
$F: k^2 \rightarrow k^2$ is tame,
meaning that it can be written as 
a finite composition of elementary maps
$(X_1 + a, X_2)$ and $(X_1,X_2+b)$.
For $n=3$ it is an open problem to determine whether
every invertible polynomial 
$F \in k[X_1, \ldots, X_n]$ is tame.
The general belief is that they are not.
\section{PERFORMANCE ISSUES}\label{sec:performanceissues}
The performance objectives of observer synthesis are:
(1) asymptotic stability of the difference of the output of the system
and the output of the observer; and
(2) a good transient response in reaction to 
realistic initial conditions of the observer.
In control theory it is known that
it is best for the functioning of the observer
if the convergence of the predicted output of the
observer to the observed output
is slightly faster than the dynamics of the system.
It is well known from other subareas of control theory
that for the analysis of the performance of an observer
one has to consider the performance system with as state $(x, x_o)$.
\par
The approach to investigate the global asymptotic stability
of the observer includes the steps:
(1) Investigate local stability at zero error.
(2) Approximate the domain of attraction by analysis
and simulation.
See \cite{vannelli:vidyasagar:1985}
for a procedure to approximate the domain of attraction
by rational Lyapunov functions.
\begin{definition}\label{def:performancesystem}
Define the {\em performance system} of the observer
as the control system,
\begin{eqnarray*}
     x_e(t)
  & = & \left(
        \begin{array}{l}
          x(t) \\ x_o(t)
        \end{array}
        \right), ~~ 
        x_e(0) 
        = 
        \left(
        \begin{array}{l}
          x_0 \\
          x_{o,0}
        \end{array}
        \right), \\
	\frac{dx_e(t)}{dt}
  & = & f_e(x_e(t))
        =
        \left(
        \begin{array}{l}
          f(x(t)) \\ f_o(x_o(t), C x(t))
        \end{array}
        \right), \\
      e_y(t)
  & = & y(t) - h_o(x_o(t)) = h(x(t)) - h_o(x_o(t)) \\
  & = & h_e(x_e(t)).
\end{eqnarray*}
\end{definition}
Note that the dimensions of the vectors
$x$ and $x_o$ are in general different
hence it is mathematically not possible to subtract these vectors.
\begin{problem}
Can an observer be synthesized such that
for any initial condition $x_e(0)$
one has that,
\[
  \lim_{t \rightarrow \infty} e_y(t) = 0? 
\]
Equivalently,
is the observable part of the performance system 
such that the output $e_y$ is asymptotically stable 
for all initial conditions of the performance system?
\end{problem}
The problem above is briefly analyzed for linear systems.
\begin{definition}
Performance system of a linear observer.
\begin{eqnarray*}
      \frac{dx(t)}{dt}
  & = & A x(t), ~ x(0) = x_0,  ~~
        y(t) = C x(t).\\
      \frac{dx_o(t)}{dt}
  & = & A x_o(t) + K [y(t) - C x_o(t)], ~ x_o(0) = x_{o,0}.\\
     x_e(t)
  & = & \left(
        \begin{array}{ll}
          x(t) & x_o(t)
        \end{array}
        \right)^T, \\
      \frac{dx_e(t)}{dt}
  & = & \left(
        \begin{array}{ll}
          A & 0 \\
          - K C & A - K C
        \end{array}
        \right) x_e(t), \\
      e_y(t)
  & = & y(t) - C x_o(t), \\
      x_{et}(t)
  & = & \left(
        \begin{array}{ll}
          x(t) & x(t) - x_o(t)
        \end{array}
        \right)^T
        = L_e x_e(t), \\
      \frac{dx_{et}(t)}{dt}
  & = & \left(
        \begin{array}{ll}
          A      & 0 \\
          0      & A - K C
        \end{array}
        \right) x_{et}(t), \\
      e_y(t)
  & = & y(t) - C x_o(t) =
        \left(
        \begin{array}{ll}
          0 & C 
        \end{array}
        \right) x_{et}(t).
\end{eqnarray*}
As is well known in control theory,
if the tuple $(A, C)$ is an observable pair
then there exists a gain matrix $K$ such that
the observable part of the performance system
is globally asymptotically stable.
Note the nonobservability of the performance system with respect to $e_y$!
\end{definition}
\begin{example}\label{ex:polynomialsyspolynomialobs}
This is a continuation of Example \ref{ex:polsyspolobs5}.
If $k_1 < - a_{22}$ and if $k_2 < 0$
then the observable part of the performance system 
is locally stable at $x - x_o = 0$. 
The system itself is such that its linearized system at
$x=0$ has one eigenvalue at zero.
The global asymptotic stability is not yet established.
\end{example}
\section{EXAMPLES}\label{sec:examples}
\begin{example}\label{ex:observerhigherorder}
{\em A polynomial system with a polynomial observer of higher
state-space dimension than the system}. 
There follows an example of an observer for a system
in which the state-space dimension of the observer
is strictly larger than that of the system.
\par
Consider the polynomial system and derive
its observer according to the steps,
\begin{eqnarray*}
     \frac{dx(t)}{dt}
  & = & \left(
        \begin{array}{l}
          f_1(x(t)) \\
          - a_{21} x_2(t) + a_{22}
        \end{array}
        \right), \\
  &   & x(0) = x_{0}, \\
      y(t)
  & = & x_1(t) = C x(t) =
        \left(
        \begin{array}{ll}
          1 & 0
        \end{array}
        \right) x(t), \\
      f_1(x) 
  & = & 2 a_{21} x_1 - a_{12} (x_2 - a_{13})(x_2 - a_{14}), \\
  &   & a_{12}, ~ a_{13}, ~ a_{14}, ~ a_{21}, ~ a_{22}  \in (0,\infty), \\
  &   & a_{13} < \frac{a_{22}}{a_{21}} < a_{14}, ~
        a_{13} + a_{14} \neq 2 \frac{a_{22}}{a_{21}}. \\
      s_1(x)
  & = & h(x) = C x = x_1, \\
      s_2(x)
  & = & \frac{\partial h(x)}{\partial x} f(x) = f_1(x), \\
     \hat{x}_1
  & = & s_1(x) = x_1, ~~
        x_1 = \hat{x}_1, \\
      \hat{x}_2
  & = & s_2(x) = 2 a_{21} x_1 - a_{12} (x_2 - a_{13})(x_2 - a_{14}) \\
  & = & 2 a_{21} \hat{x}_1 - a_{12} (x_2 - a_{13})(x_2 - a_{14}).
\end{eqnarray*}
The conclusion is that based on the state-space transformation
$s = ( s_1, ~ s_2)$
there does not exist an unique solution for the inverse function,
in particular not a solution for $x_2$ given $\hat{x}_1, ~ \hat{x}_2$
of the equation
\[
  a_{12} (x_2 - a_{13}) (x_2 - a_{14}) = - \hat{x}_2 + 2 a_{21} \hat{x}_1.
\]
Therefore the construction of the state-space transformation
is increased by another derivation.
\begin{eqnarray*}
      s_3(x)
  & = & \frac{\partial s_2(x)}{\partial x} f(x) 
        = c_{11} x_1 + c_{12} x_2 + c_{14}, \\
      c_{11}
  & = & 4 a_{21}^2, \\
      c_{12}
  & = & a_{12} a_{21} (a_{13} + a_{14}) 
        - 2 a_{12} a_{22}, \\
      c_{14}
  & = & a_{12} a_{22} (a_{13} + a_{14}) 
        - 2 a_{12} a_{21} a_{13} a_{14}. \\
      \hat{x}_3
  & = & s_3(x), \\
      x_2
  & = & - \frac{c_{11}}{c_{12}} \hat{x}_1
        + \frac{1}{c_{12}} \hat{x}_3 
        - \frac{c_{14}}{c_{12}}, \\
      (x_1, ~ x_2)
  & = & s^{-1}(\hat{x}) 
        = (\hat{x}_1, ~ 
             - \frac{c_{11}}{c_{12}} \hat{x}_1
             + \frac{1}{c_{12}} \hat{x}_3 
             - \frac{c_{14}}{c_{12}}
          ), ~\\
      \frac{d \hat{x}(t)}{dt}
  & = & A_{or} \hat{x}(t) + b_{or}(\hat{x}(t)), ~~
        \hat{x}(t), ~ x_o(t) \in \mathbb{R}^3.
\end{eqnarray*}
Note that $s^{-1}$ is a linear map. 
Then,
\begin{eqnarray*}
      s_4(x)
  & = & 2 a_{21} c_{11} x_1
        - a_{21} c_{12} x_2 + \\
  &   & - a_{12} c_{11} (x_2 - a_{13}) ( x_2 - a_{14})
        + a_{22} c_{12},  \\
      b_{or,3}(\hat{x})
  & = & s_4(x)|_{x = s^{-1}(\hat{x})} \\
  & = & 3 a_{21} c_{11} \hat{x}_1 
        - a_{21} \hat{x}_3
        + a_{21} c_{14} 
        + a_{22} c_{12} +  \\
  &   & - \frac{a_{12} c_{11}}{c_{12}^2} \times \\
  &   & \times 
          [ \hat{x}_3 - c_{11} \hat{x}_1 - ( c_{14} + a_{13} c_{12}) ] \times \\
  &   & \times
          [ \hat{x}_3 - c_{11} \hat{x}_1 - ( c_{14} + a_{14} c_{12}) ], 
\end{eqnarray*}
\begin{eqnarray*}
      k_{o,3}(\hat{x})
  & = & \frac{\partial b_{or,3}(\hat{x})}{\partial \hat{x}_1} 
        = 3 a_{21} c_{11} 
        + \frac{a_{12} c_{11}^2}{c_{12}^2} \times \\
  &   & \times
        \left[
        \begin{array}{l}
          + 2 \hat{x}_3
          - 2 c_{11} \hat{x}_1 + \\
          - ( 2 c_{14} + c_{12} ( a_{13} + a_{14} ) )
        \end{array}
        \right].
\end{eqnarray*}
The observer is then described by the system representation,
\begin{eqnarray*}
     \frac{dx_o(t)}{dt}
  & = & A_{o} x_o(t) + b_{o}(x_o(t))  + \\
  &   & + [ k_o(x_o(t)) + K ] [ y(t) - C_o x_o(t) ] \\
  & = & \left(
        \begin{array}{l}
          \hat{x}_2(t) \\
          \hat{x}_3(t) \\
          b_{o,3}(\hat{x}(t))
        \end{array}
        \right) + \\
  &   & +
        \left(
        \begin{array}{l}
          k_1 \\
          k_2 \\
          k_3 + k_{o,3}(\hat{x}(t))
        \end{array}
        \right) [ y(t) - C_o \hat{x}(t) ],  \\
      y_o(t)
  & = & C_o x_o(t), ~~
     C_o
     =  \left(
        \begin{array}{lll}
          1 & 0 & 0
        \end{array}
        \right). 
\end{eqnarray*}
\end{example}
\begin{example}\label{ex:ratsysobsratsys}
{\em Rational system -- rational observer}.
Consider the observable rational system,
for $a_{13} \neq 0$.
\begin{eqnarray*}
     \frac{dx(t)}{dt}
  & = & \left(
        \begin{array}{l}
          \frac{-a_{11}x_1(t)}{1+a_{12}x_1(t)}
          + \frac{a_{13}x_2(t)}{1+a_{14}x_2(t)} \\
          \frac{-a_{21}x_2(t)}{1+a_{22}x_2(t)} + a_{23}
        \end{array}
        \right), \\
      y(t)
 & = & C x(t) =
       \left(
       \begin{array}{ll}
         1 & 0 
       \end{array}
       \right) x(t) = x_1(t).
\end{eqnarray*}
Construct an output-based realization according
to Procedure \ref{proc:observersynthesis}.
\begin{eqnarray*}
      \frac{d \hat{x}(t)}{dt} 
  & = & A_{or} \hat{x}(t) 
        + b_{or}(\hat{x}(t))  \\
  & = & \left(
        \begin{array}{ll}
          0 & 1 \\
          0 & 0
        \end{array}
        \right) \hat{x}(t)
        +
        \left(
        \begin{array}{l}
          0 \\ b_{or,2}(\hat{x}(t))
        \end{array}
        \right), \\
      y(t)
  & = & C \hat{x}(t), \\
      b_{or,2}(x)
  & = & \frac{a_{11}^2 x_1}{(1+a_{12}x_1)^2} 
        - \frac{a_{11} a_{13} x_2}{(1+a_{12}x_1)^2(1+a_{14}x_2)} + \\
  &   & - \frac{a_{13}a_{21} x_2}{(1+a_{14}x_2)^2 (1+a_{22}x_2)}
        + \frac{a_{13}a_{23}}{(1+a_{14}x_2)^2}. 
\end{eqnarray*}
The rational observer follows.
\begin{eqnarray*}
  \frac{dx_o(t)}{dt}
  & = & A_o x_o(t)
        + b_o (x_o(t)) + \\
  &   & + [ k_o(x_o(t)) + K ] [ y(t) - C x_o(t)], \\
      A_o
  & = & A_{or}, \\
      b_o(x_o)
  & = & b_{or}(x_o) =
        \left(
        \begin{array}{l}
          0 \\ b_{or,2}(x_o)
        \end{array}
        \right), \\
      \lefteqn{
        k_{o,2}(x_o)
        = \frac{\partial b_{or,2}(x_o)}{\partial x_{o,1}} 
      } \\
  & = & \frac{a_{11}^2}{(1+a_{12}x_{o,1})^3} 
        - \frac{3 a_{11} a_{12} x_{o,1}}{(1 + a_{12} x_{o,1})^4} + \\
  &   & + \frac{2 a_{12}}{(1 + a_{12}x_{o,1})^3} 
        \frac{a_{11} a_{13} x_{o,2}}{(1 + a_{14} x_{o,2})} + \\
  &   & - \frac{a_{13} a_{14} x_{o,2}}{(1 + a_{14} x_{o,2})^2 (1 + a_{22} x_{o,2})}.
\end{eqnarray*}
\end{example}
\begin{example}\label{ex:brs2}
Consider a rational system which is a model
for either a compartmental system with two compartments
or a biochemical reaction system with two chemical species.
The output is only of the second compartment or species.
The system is actually a family of rational systems.
\begin{eqnarray*}
      \frac{dx_1(t)}{dt}
  & = & r_1(x_1(t)) - r_2(x_1(t),x_2(t)), ~ x_1(0) = x_{1,0},\\
      \frac{dx_2(t)}{dt}
  & = & r_2(x_1(t),x_2(t)) - r_3(x_2(t)), ~ x_2(0) = x_{2,0}, \\
      y(t)
  & = & C x(t) = x_2(t), ~ 
        r_1, ~ r_2, ~ r_3 \in Q_X.
\end{eqnarray*}  
Assume that
$r_2(0,x_2) = 0$ and $r_3(0) =0$ for all $x_2 \in [0,\infty)$
and that $x_0 = (x_{0,1}, x_{0,2}) \in \mathbb{R}_+^2$.
From those assumptions follows that
the positive orthant $\mathbb{R}_+^2 = [0,\infty)^2$
is an invariant set 
to be denoted as the state set $X = \mathbb{R}_+^2$.
A special case of the above is with
\begin{eqnarray*}
  &   & r_1(x_1) = a_1, ~
        r_2(x_1,x_2) = \frac{a_2 x_1 x_2}{1+ a_3 x_2^2}, \\
  &   & r_3(x_2) = a_4 x_2, \quad 
        a_1, a_2, a_3, a_4 \in (0,\infty).
\end{eqnarray*}
The state-space transformation and its inverse are then
\begin{eqnarray*}
      \hat{x}_1
  & = & s_1(x) = h(x) = C x = x_2, \\
      \hat{x}_2
  & = & s_2(x) = \frac{\partial s_1(x)}{\partial x}
        = f_2(x) = r_2(x_1,x_2) - r_3(x),\\
      x_2
  & = & \hat{x}_1, \quad 
        x_1 = \frac{(\hat{x}_2 - a_4 \hat{x}_1) (1+a_3 \hat{x}_1^2)}{a_2 \hat{x}_1}.
\end{eqnarray*}
The output-based realization is then,
\begin{eqnarray*}
     \frac{d \hat{x}(t)}{dt}
  & = & \left(
        \begin{array}{l}
          \hat{x}_2(t) \\ b_2(\hat{x}(t))
        \end{array}
        \right), ~ 
        \hat{y}(t) = C \hat{x}(t), \\
      b_{o,2}(\hat{x})
  & = & p_{o,2}(x_o)/[ a_2 x_{o,1} (1+a_3 x_{o,1}^2) ], \\
      p_{o,2}(x_o)
  & = & - 2 a_2 x_{o,1}^2 + a_4^2 x_{o,1} ( 1 + a_3 x_{o,1}^2 )^2 + \\
  &   & + (1+a_3x_{o,1}^2)  
        [ a_2 x_{o,1} ( a_1 - 2 a_4 x_{o,1}) + \\
  &   & + 2 (x_{o,2} + a_4 x_{o,1} ) + \\
  &   & -2 a_1 x_{o,1}^4 ( x_{o,2} + a_4 x_{o,1})^2  
        ].
\end{eqnarray*}
The observer is then,
\begin{eqnarray*}
      \frac{dx_o(t)}{dt}
  & = & \left(
        \begin{array}{l}
          x_{o,2}(t) \\ b_o(x_o(t))
        \end{array}
        \right) + \\
  &   & + [ k_o(x_o(t)) + K ] ( y(t) - C_o x_o(t)).
\end{eqnarray*}
\end{example}
\section{CONCLUDING REMARKS}\label{sec:concludingremarks}
The main contribution of the paper is 
an observer synthesis procedure for rational observers
of rational systems.
\par
Further research into the subject of this paper includes:
Synthesis of observers for rational systems with inputs 
(first for systems affine in the input).
The role of universal inputs has to be explored.
A stability analysis of rational observers.
Observers for Nash systems 
\cite{nemcova:petreczky:schuppen:2013:siamjco}
have been constructed but are not
discussed in this paper due to space limitations.




%
\begin{footnotesize}
%

%
\end{footnotesize}

\end{document}